\documentclass{amsart}
\usepackage{bbm}
\usepackage{graphicx}
\usepackage{verbatim}
\usepackage{stackrel}
\usepackage{color}
\graphicspath{{Fig/}}
\usepackage{amsmath}
\usepackage{mathtools}
\usepackage{enumerate}

\usepackage[latin1]{inputenc}
\def\diag{\text{\rm Diag}}

\newtheorem{theorem}{Theorem}

\theoremstyle{definition}
\newtheorem{prop}{Proposition}
\newtheorem{lemma}{Lemma}

\newtheorem{defi}{Definition}

\newtheorem*{assumption}{Assumptions}

\def\N{{\mathbb N}}
\def\R{{\mathbb R}}
\def\Z{{\mathbb Z}}
\def\P{{\mathbb P}}
\def\E{{\mathbb E}}
\def\Z{{\mathbb Z}}

\newcommand{\diff}{\mathop{}\mathopen{}\mathrm{d}}

\newcommand\ind[1]{\mathbbm{1}_{\left\{#1\right\}}}

\renewcommand\vec[1]{#1}

\def\cal{\mathcal}

\title{A functional central limit theorem for the Becker-D\"oring model}

\address[Wen SUN]{INRIA Paris, 2 rue Simone Iff, F-75012 Paris, France}

\author{Wen Sun}
\email{Wen.Sun@inria.fr}

\date{\today}

\begin{document}

\begin{abstract}
  We investigate the fluctuations of the stochastic Becker-D\"oring model of polymerization when the initial size of the system converges to infinity.  A  functional central limit problem  is proved for the vector of the number of polymers of a given size.  It is shown that the stochastic process associated to fluctuations is converging to the strong solution of an infinite dimensional stochastic differential equation (SDE)  in a  Hilbert space.
 We also prove that, at equilibrium,  the solution of this SDE  is a  Gaussian process. The proofs are based on a  specific representation of the evolution equations,  the introduction of a convenient Hilbert space and several technical estimates to control the fluctuations, especially of the first coordinate which interacts with all components of the infinite dimensional vector representing the state of the process. 
\end{abstract}

\maketitle

\vspace{-5mm}

\bigskip

\hrule

\vspace{-3mm}

\tableofcontents

\vspace{-1cm}

\hrule

\bigskip

\section{Introduction}
Polymerization is a key phenomenon in several important biological processes. Macro-molecules proteins, also called monomers, may be assembled randomly into several aggregated states called polymers or clusters. These clusters can themselves be fragmented into monomers and polymers at some random instants. The fluctuations of the number of polymerized monomers analyzed in this paper is an important characteristic of polymerization processes in general. 

\subsection*{The  Becker-D\"oring model} 
We investigate the fluctuations of  a stochastic version of the Becker-D\"oring model which is a classical mathematical model to study polymerization.  The Becker-D\"oring model describes the time evolution of the distribution of cluster sizes in a system where only additions (coagulation mechanism) or removals (fragmentation) of one monomer from a cluster are possible.  A cluster of size~$1$ is a monomer and  clusters  of size greater than $2$ are  polymers.  Under Becker-D\"oring model, coagulation and fragmentation are  simple synthesis and decomposition reactions: a polymer of size $k$ may react with a monomer to form a polymer of size $k{+}1$ at kinetic rate $a_k$; a polymer of size $k{+}1$ may break down into a polymer of size $k$ and a monomer at kinetic rate $b_{k{+}1}$, i.e.,
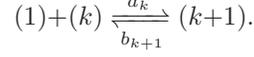
\[
(1){+}(k)\xrightleftharpoons[b_{k{+}1}]{\,a_k\,}(k{+}1).
\]
The ODEs associated to the deterministic version of the Becker-D\"oring model have been widely studied in physics since 1935, see Becker and D\"oring~\cite{BDeq1935}. This is an infinite system of ordinary differential equations of $c(t){=}(c_k(t),k{\in}\N^+)$, given by 
 \begin{equation}\tag{BD}
   \begin{cases}
\displaystyle     \dfrac{  \diff c_1}{\diff t} (t){=}{-}2J_{1}(c(t)){-}\sum_{k\geq2} J_{k}(c(t)), \\ \dfrac{  \diff c_k}{\diff t} (t){=}J_{k-1}(c(t)){-}J_{k}(c(t)), \, k{>}1,
   \end{cases}
\end{equation}
 with  $J_k(c){=}a_kc_1c_k{-}b_{k+1}c_{k+1}$ if $c{=}(c_k){\in}\R_+^{\N^+}$.  For $k{\geq}1$, $c_k(t)$ represents the concentration of clusters of size~$k$ at time $t$. 
 The conditions on existence/uniqueness of solutions  for the Becker-D\"oring equations (BD) have been extensively investigated. See Ball et al.~\cite{Ball1986},  Niethammer~\cite{niethammer2004} and  Penrose~\cite{Penrose,Penrose2}.
 
The evolution equations satisfied by  $(c(t))$ can be rewritten under a more compact form as,
\begin{align}\label{bdeq}
\diff c(t){=}\tau{\circ} s(c(t))\diff t,
\end{align}
where $s$ is a mapping from $\R_+^{\N^+}{\to}\R_+^{\N^+}$: for any $k{\in}\N^+$ and $c{\in} \R_+^{\N^+}$
\begin{equation}\label{defs}
s_{2k-1}(c){=}a_k c_1c_k\quad\textrm{ and  }\quad
s_{2k}(c){=}b_{k+1} c_{k+1};
\end{equation}
and $\tau$ is a linear mapping from $\R_+^{\N^+}{\to}\R_+^{\N^+}$: for any $z{\in} \R_+^{\N^+}$ and $k{\ge} 2$,
\begin{equation}\label{deftau}
\begin{cases}
\displaystyle  \tau_1(z){=}{-}\sum_{i{\ge} 1}(1{+}\ind{i=1})z_{2i{-}1}{+}\sum_{i{\ge} 2}(1{+}\ind{i{=}2})z_{2i{-}2},\\
\tau_k(z){=} z_{2k{-}3}{-}z_{2k{-}2}{-}z_{2k{-}1}{+}z_{2k}.
\end{cases}
\end{equation}
As it will be seen this representation will  turn out to be very useful to derive the main results  concerning fluctuations.

\subsection*{Becker-D\"oring ODEs and  Polymerization Processes}
This set of ODEs is  used to describe the evolution of the concentration $c_i(t)$, $i{\ge}1$,  of polymers of size~$i$.  The classical framework assumes an initial state with only polymers of size $1$, monomers. In a biological context, experiments show that the concentration of polymers of size greater than $2$ stays at $0$ until some instant, defined as the {\em lag time }, when the polymerized  mass grows very quickly to reach its stationary value.  With convenient parameters estimations, these ODEs can be used to describe first order characteristics such as mean concentration of polymers of a given size.  The use of systems of ODEs to describe the evolution of polymerization processes started with  Oosawa's pioneering work  in 1962, see Oosawa and Asakura~\cite{Oosawa} for example.  Morris et al.~\cite{Morris} presents  a quite detailed review of the classical sets of ODEs used for polymerization processes. As it can be seen and also expected, the basic dynamics of the Becker-D\"oring model of adding/removing a monomer to/from a polymer occupy a central role in most of these mathematical models. See also Prigent et al.~\cite{Prion} and Hingant and Yvinec~\cite{Yvinec2016} for recent developments in this domain. 

Outside the rapid growth of polymerized mass at the lag time, the other important aspect observed in the experiments is the high variability of the instant when it occurs, the lag time, from an experiment to another.  This is believed to explain, partially,  the variability of the starting point of  diseases associated to these phenomena in neural cells, like Alzheimer's disease for example.  See Xue et al.~\cite{Radford}. Hence if the deterministic Becker-D\"oring ODEs describes the first order of polymerization through a limiting curve, the fluctuations around these solutions  will give a characterization of the variability of the processes itself. Up to now the mathematical studies of these fluctuations are quite scarce, the stochastic models analyzed include generally  only a finite number of possible sizes for the polymers.    See  Szavits et al.~\cite{Szavits}, Xue et al.~\cite{Radford} and Eug\`ene et al.~\cite{EXRD} for example.  To study these important aspects,  we have to introduce a stochastic version of the Becker-D\"oring model. 

\subsection*{The stochastic Becker-D\"oring model}
The polymerization process is described as a Markov process $(X^N(t)){:=}(X_k^N(t), 1 {\le} k{\le} N)$,
where $X_k^N(t)$ is  the total number of clusters of size $k$  at time $t$,   it takes values in the state space,
    \[
  \cal{S}^N{:=}\bigg\{\vec{x}{\in}\N^N\bigg| \sum_{k= 1}^Nkx_k{=} N\bigg\},
  \]
of configurations of polymers  with  mass $N$.

For any $\vec{x}{\in}\cal{S}^N$ and $1{\le} k{<}N$, the associated jump matrix $Q^N{=}(q^N(\cdot,\cdot))$  is given by
\begin{displaymath}
\begin{cases}
 q^N(\vec{x},\vec{x}{-}\vec{e}_1{-}\vec{e}_k{+}\vec{e}_{k{+}1}){=}a_k x_1(x_k{-}\ind{k{=}1})/N, \\
 q^N(\vec{x},\vec{x}{+}\vec{e}_{1}{+}\vec{e}_k{-}\vec{e}_{k{+}1}){=}b_{k{+}1}x_{k{+}1}, 
\end{cases}
\end{displaymath}
where $(\vec{e}_k,k{\in}\N^+)$ is the standard orthonormal basis of $\N^{\N^+}$. In other words, a monomer is added to a given polymer of size $k$ at rate $a_kx_1/N$ and a monomer is detached from a polymer of size $k$ at rate $b_k$.  Note that if $N$ is interpreted as a ``volume'', the quantity $X_1^N(t)/N$ can be seen as the concentration of monomers. 

There are few studies of this important stochastic process.  The large scale behavior of the stochastic Becker-D\"oring model, when $N$ gets large, is an interesting and challenging problem. 
Given the transition rates, one can expect that the deterministic Becker-D\"oring equations~(BD)   give the limiting equations for the concentration of the different species of polymers, i.e. for the convergence in distribution
\[
\lim_{N\to+\infty}(X_k^N(t)/N, k{\ge}1){=}(c_k(t),k\geq1).
\]
Such a first order analysis is achieved in Jeon~\cite{Jeon1998}. This result is in fact proved for a more general model, the Smoluchowski coagulation-fragmentation model.
The Becker-D\"oring model is a special case of Smoluchowski model, see  Aldous~\cite{Aldous1999} for a survey on the coalescence models.
The stochastic approximation of the pure Smoluchowski coalescence equation through Marcus-Lushnikov process is investigated in Norris~\cite{norris1999}.
There are also results on the relation of stochastic Becker-D\"oring model and the deterministic Lifshitz-Slyozov equation, see~Deschamps et al.~\cite{deschamps} for example.

\subsection*{The Main Contributions}
In this paper,  we investigate the fluctuations of the Becker-D\"oring model, i.e.  the limiting behavior of the $\R^N$-valued process
\begin{equation}\label{Wn}
\left(W^N(t)\right){:=}\left(\frac{1}{\sqrt{N}}\left(X^N(t){-}Nc(t)\right)\right)
\end{equation}
is analyzed. We prove that, under appropriate conditions,  the fluctuation process $(W^N(t))$ converges for the Skorohod topology  to a ${L_2(w)}$-valued process $(W(t))$  which is the  strong solution of the  SDE
\begin{equation}\label{limitsde1}
 \diff W(t){=}\tau\bigg(\nabla s(c(t))\cdot W(t)\bigg)\diff t{+} \tau\bigg(\diag\left(\sqrt{s(c(t))}\right)\cdot\diff \beta (t)\bigg),
\end{equation}
where
\begin{enumerate}[{\qquad 1.}]
  \item ${L_2(w)}$ is a Hilbert subspace  of $\R^{\N^+}$, see definition~\ref{defH} of Section~\ref{limit};
\item The operators $s(\cdot)$ and  $\tau(\cdot)$  are defined respectively by relations~\eqref{defs} and~\eqref{deftau};
\item $\nabla s(c)$ is the Jacobian matrix
  \[
\nabla s(c){=}\left(\frac{\partial s_i}{\partial c_j}(c)\right), c{\in}\R_+^{\N^+} 
\]
\item $\beta(t){=}(\beta_k(t),k{\in}\N^+)$ is a sequence of independent standard Brownian motions in $\R$;
\item $\diag(\vec{v})$ is a diagonal matrix whose diagonal entries  are the coordinates of the vector $\vec{v}{\in} \R^{\N^+}$.
\end{enumerate}
See Theorem~\ref{fclt} in Section~\ref{limit} for a precise formulation of this result.  Note that  the drift part of~\eqref{limitsde1} is the gradient of the Becker-D\"oring equation~\eqref{bdeq}.

Under appropriate assumptions, see Ball et al.~\cite{Ball1986}, the  Becker-D\"oring equations \eqref{bdeq} have a unique fixed point  $\widetilde{c}$. 
If $c(0){=}\widetilde{c}$, the asymptotic  fluctuations around this stationary state are then given  by $(\widetilde{W}(t))$. It is shown in Proposition~\ref{fe} that this process, the solution of SDE~\eqref{limitsde1},  can be represented as 
\[
 \widetilde{W}(t){=}\cal{T}(t)\widetilde{W}(0){+}\int_0^t\cal{T}(t{-}s)\tau\Big(\diff \cal{B}(s)\Big),
\]
  where $\cal{T}$ is the  semi-group associated with linear operator $\tau{\circ}\nabla s(c)$ and $(\cal{B}(t))$ is a $\widetilde{Q}$-Wiener process in ${L_2(w)}$ where
  $\widetilde{Q}{=} \diag\left(w_ns_n(c),n{\ge}1\right)$. In particular if the initial state $\widetilde{W}(0)$ is deterministic, then $( \widetilde{W}(t))$ is a Gaussian process.
  For the definitions and properties of $\widetilde{Q}$-Wiener process and Gaussian processes in Hilbert spaces, see Section~4.1 and Section~3.6 in Da Prato and Zabczyk~\cite{Prato1992} for example.

  \subsection*{Literature}
For the fluctuation problems in the models with finite chemical reactions, results are well known, see Kurtz~\cite{kurtz1,kurtz2} for example. However, in the Becker-D\"oring model, there are countable many species and reactions, where the results for the finite reactions are not directly applicable. See the Open Problem~9 in Aldous~\cite{Aldous1999} for the description of fluctuation problem in the general Smoluchowski coagulation model.
  For this reason, the studies of fluctuations of the Becker-D\"oring models are quite scarce. Ranjbar and Rezakhanlou~\cite{Ranjbar2010} investigated the fluctuations at equilibrium  of a family of coagulation-fragmentation models  with a spatial component. They proved that, at equilibrium,  the limiting fluctuations can be described  by an Ornstein-Uhlenbeck process in some abstract space. Their results rely on balance equations which hold because of  the stationary framework.
  Durrett et al.~\cite{durrett} gave the stationary distributions for all reversible coagulation-fragmentation processes. Then they provided the limits of the mean values, variances and covariances of the stationary densities of particles of any given sizes when the total mass tends to infinity.

  Concerning fluctuations of infinite dimensional Markov processes, several examples have received some attention, in statistical physics mainly.
In the classical Vlasov model, the first result on the central limit theorem seems to be given by Braun and Hepp~\cite{braun} in 1977.
They investigated the fluctuations around the trajectory of a test particle.
The central limit theorem for the general McKean-Vlasov model when the initial measures of the system are products of i.i.d. measures is proved by Sznitman~\cite{SZNITMANclt} in 1984. For the Ginzburg-Landau model on $\Z$, where the evolution of the state at a site $i$  depends only on the state of its nearest neighbors $i{\pm}1$.  The independent (Brownian) stochastic fluctuations at each site do not depend of the state of the process.  The hydrodynamic limit, i.e. a first order limit,  is given by Guo, Papanicolao and Varadhan~\cite{guo1988}.  The fluctuations around this hydrodynamic limit live in an infinite dimensional space.   Zhu~\cite{zhu1990} proved that, at equilibrium,  the limiting fluctuations converge to  a stationary Ornstein-Uhlenbeck process taking values in the dual space of a nuclear space.  In  the non-equilibrium case, the fluctuations have been investigated  by Chang and Yau~\cite{chang1992},  the limiting fluctuations can be described as an Ornstein-Uhlenbeck in a negative Sobolev space.  For more results on related fluctuation problems in statistical mechanics, see Spohn~\cite{spohn1986,spohn2012}.  Fluctuations of an infinite dimensional Markov process associated to  load balancing mechanisms in large stochastic networks have been investigated in Graham~\cite{Graham2005}, Budhiraja and Friedlander~\cite{budhiraja2017}.

One of the difficulties of our model is the fact that the first coordinate $(X_1^N(t))$ of our Markov process,  the number of monomers,   interact with   all non-null coordinates of the state process. Recall that monomers can react with all other kinds of polymers. This feature has implications on the choice of the Hilbert space ${L_2(w)}$ chosen to formulate the SDE~\eqref{limitsde} and, additionally,  several estimates have to be derived to control the stochastic fluctuations of the first coordinate. This situation is different from the examples of the Ginzburg-Landau model in~\cite{chang1992,spohn1986,zhu1990}, since  each site only have interactions with a finite number of sites, or in the stochastic network example~\cite{budhiraja2017} where the interaction range is also finite. It should be also noted that our evolution equations are driven by a set of independent Poisson processes whose intensity is state dependent which is not the case in the Ginzburg-Landau models for which the diffusion coefficients are constant.  For this reason Lipschitz properties have to be established for an appropriate norm for several functionals, it complicates the already quite technical framework of these problems. 

\subsection*{Outline of the paper}
Section~\ref{model} introduces the notations, assumptions and the stochastic model as well as the evolution equations.  Section~\ref{limit} investigates the problem of existence and uniqueness of the solution of the SDE~\eqref{limitsde1}.  The proof of the fluctuations at equilibrium is also given.  Section~\ref{convergence} gives the proof of the main result, Theorem~\ref{fclt}, the convergence of the fluctuation processes.

\section{The Stochastic Model}\label{model}
In this section we introduce the notations and assumptions used throughout this paper.
The stochastic differential equations describing the evolution of the model are introduced.

For any $h{\in}\R^+$,  $(\cal{N}_h^i)_{i=1}^\infty$ denotes  a sequence of independent Poisson processes with intensity $h$. The stochastic Becker-D\"oring equation  with aggregation rates $(a_k,k\in\N^+)$ and fragmentation rates $(b_k,k\in\N^+)$ can be expressed  as the solution of the SDEs
  \begin{align*}
    \diff X_1^N(t){=}&{-}\sum_{k{\ge} 1}(1{+}\ind{k=1})\hspace{-8mm}\sum_{i=1}^{X_1^N(t-)(X_k^N(t-){-}\ind{k{=}1})}\hspace{-8mm}\cal{N}_{a_k/N}^i(\diff t)
    +\sum_{k\ge 2}(1{+}\ind{k=2})\sum_{i=1}^{X_k^N(t-)}\!\!\!\cal{N}_{b_k}^i(\diff t),\\
    \diff X_k^N(t){=}&\hspace{-6mm}\sum_{i{=}1}^{X_1^N(t-)(X_{k-1}^N(t-){-}\ind{k{=}2})}\hspace{-8mm}\cal{N}_{a_{k-1}/N}^i(\diff t){-}\sum_{i=1}^{X_k^N(t-)}\cal{N}^i_{b_k}(\diff t)\\
    &\qquad {-}\hspace{-4mm}\sum_{i=1}^{X_1^N(t-)X_{k}^N(t-)}\hspace{-4mm}\cal{N}_{a_{k}/N}^i(\diff t){+}\sum_{i=1}^{X_{k+1}^N(t-)}\cal{N}^i_{b_{k+1}}(\diff t),
  \end{align*}
and for all $k>N$, $X_k^N(t){\equiv} 0$, with $f(t{-})$ being the limit on the left of the function $f$ at $t{>}0$.

 In order to separate the drift part and the martingale part, it is convenient to introduce the corresponding martingales $(D^N(t))$.
  For $k{\ge}1$, let
\begin{align*}
  & D_{2k-1}^N(t)=\int_0^t\left(\sum_{i=1}^{X_1^N(u-)(X_k^N(u-){-}\ind{k{=}1})\hspace{-10mm}}\cal{N}^i_{a_k/N}(\diff u)-\frac{1}{N}a_kX_1^N(u)(X_k^N(u){-}\ind{k{=}1})\diff u\right),\\
    &D_{2k}^N(t)=\int_0^t\left(\sum_{i=1}^{X_{k+1}^N(u-)}\cal{N}^i_{b_{k+1}}(\diff u)-b_{k+1}X_{k+1}^N(u)\diff u\right).\\
\end{align*}
Clearly, for every fixed $i{\in}\N^+$, $(D_{i}^N(t))$  is local martingale in $\R$ with previsible increasing process
\begin{equation}\label{inc}
\left(\langle D_{i}^N\rangle (t)\right)=\left(N\int_0^ts_{i}\left(\frac{X^N(u)}{N}\right)\diff u-\ind{i=1}\int_0^ta_1\frac{X_1^N(u)}{N}\diff u\right),
\end{equation}
where the operator $s=(s_i, i{\ge}1)$ is defined by relation~\eqref{defs}, additionally the cross-variation processes are null, i.e., for any $i{\neq} j$,
\[
\left(\langle D_i^N,D_j^N \rangle (t)\right)=(0). 
\]
This is in fact one of the motivations to introduce the variables $(D_k^N(t))$ and the functionals $(\tau(\cdot))$ and $(s(\cdot))$.

A simple calculation shows that, for any $N$, the process $(X^N(t))$ satisfies the relation
\begin{multline}\label{bdsde}
  X^N(t){=}  X^N(0){+}N\int_0^t\tau\left(s\left(\frac{X^N(u)}{N}\right)\right)\diff u\\{+}\tau\left(D^N(t)\right){+}2a_1\int_0^t\left(\frac{X_1^N(u)}{N}\right)e_1\diff u,
\end{multline}
where $\tau$ is defined by relation~\eqref{deftau} and $e_1{=}(1,0,0,{\dots})$.  Therefore the fluctuation process $(W^N(t))$, see relation~\eqref{Wn},  satisfies
\begin{multline}\label{fluct}
  W^N(t){=}W^N(0){+}\frac{1}{\sqrt{N}}\tau\left(D^N(t)\right){+}\frac{2a_1}{\sqrt{N}}\int_0^t\left(\frac{X_1^N(u)}{N}\right)e_1\diff u\\
 {+}\frac{1}{2}\int_0^t\tau\left(\nabla s\left(\frac{X^N(u)}{N}\right)\cdot W^N(u){+}\nabla s\left(c(u)\right)\cdot W^N(u)\right)\diff u.
\end{multline}

If we let $N$ go to infinity, then $(X^N(t)/N)$ converges to  $c(t)$ and provides that, 
for all $i\ge 1$,
\[
\lim_{N\to+\infty} \left(\langle D^N_i/\sqrt{N}\rangle(t)\right)=\int_0^ts_i(c(u))\diff u.
\]
We can expect that the process $(D^N(t)/\sqrt{N})$ converges to a stochastic integral
\[
\int_0^t\diag\left(\sqrt{s(c(u))}\right)\cdot \diff \beta(u),
\]
with $(\beta(t))$ defined in relation~\eqref{limitsde1} in  the introduction.  See page~339  in Ethier and Kurtz's~\cite{EAK}  for example. 
Then, formally, the limiting process for $(W^N(t))$, provided that it exists,  would be the solution of the SDE~\eqref{limitsde1}.

In order to give the well-posedness of the infinite dimensional process~\eqref{limitsde1}, we introduce the following notations.
\begin{defi}\label{Xplus}
  Let $\cal{X}_+$ be the phase space of the Becker-D\"oring model with initial density less than $1$,
  \[
\cal{X}_+:=\bigg\{c=\left(c_k,k\in\N^+\right)\bigg|\forall k,~c_k\ge 0;~\sum_{k\ge 1}kc_k\le 1\bigg\}.
\]
For $c(0)\in \cal{X}_+$, the solution of the Becker-D\"oring equation~(BD) is a continuous function taking values in $\cal{X}_+$.
\end{defi}
\begin{defi}\label{defH}
One assumes that $w{=}(w_n)$ is a fixed non-decreasing sequence of positive real numbers such that
  \begin{itemize}
\item  $\|1/w\|_{l_1}{:=}\sum_{n\ge 1}1/w_{n}{<}\infty$;
\item $\lim_{n\to\infty}w_{n}^{1/n}{\le} 1$;
\item there exists a constant $\gamma_0$ such that for all $n{\in}\N^+$, $w_{2n}{\le}\gamma_0 w_{n}$.

One denotes by  $L_2(w)$ the  associated $L_2$-space,
\[
L_2(w){:=}\left\{z\in\R^{\N^+}\bigg|\sum_{n=1}^\infty w_{n}z_n^2{<}\infty\right\},
\]
\end{itemize}
its inner product is defined by, for $z$, $z'{\in} {L_2(w)}$
\[
\langle z,z'\rangle_{L_2(w)}{=}\sum_{n=1}^\infty w_{n}z_nz_n',
\]
and its associated norm is $\|z\|_{L_2(w)}{=}\sqrt{\langle z,z\rangle_{L_2(w)}}$.

An  orthonormal basis $(\vec{h}_n)_{n=1}^\infty$ of $L_2(w)$ is defined as, for $n{\ge}1$,
\begin{align}\label{defhn}
  \vec{h}_n{=}(0,\dots,{1}/{\sqrt{w_{n}}},0,\dots).
\end{align}
\end{defi}
As it will be seen in Section~\ref{limit}, ${L_2(w)}$ is a convenient space to ensure a boundedness property of the linear mapping $\tau{\circ}\nabla s(\cdot)$, which is essential for the study of fluctuation process~\eqref{limitsde1}.

We now turn to the conditions on the rates and the initial state of the Becker-D\"oring equations. 
\begin{assumption}
  \begin{itemize}
    \item[]
\item[(a)]  The  kinetic rates $(a_k)$, $(b_k)$ are positive and bounded, i.e., 
\[
 \Lambda:=\max_{k\ge 1}\left\{a_k,b_k\right\}{<}\infty.
 \]
\item[(b)] There exists a positive increasing sequence  $r{\in}\R^{\N^+}_+$,   such that
  \begin{itemize}
\item  $r_{k}{\ge }w_{k}$, for all $k{\ge}1$;
\item there exists a fixed constant $\gamma_r$, such that $r_{2k}\le \gamma_r r_{k}$, for all $k{\ge}1$;
\item $(w_k/r_k)$ is a decreasing sequence  converging to $0$,
  \end{itemize}
  where $(w_n)$ is the sequence introduced in Definition~\ref{defH}.

The initial state of the  process $(X_k^N(0))$ is assumed to satisfy
\[
\sup_{N}\E\left(\sum_{k\ge 1}\frac{r_{k}X_k^N(0)}{N}\right){<}\infty.
\]
\item[(c)] The initial state of the first order process, $c(0){=}(c_k(0)){\in} \cal{X}_+$, is  such that
  \[
  \sum_{k\ge 1}w_{k}c_k(0){<}\infty \text{ and } \lim_{N\to\infty}\E\left(\sum_{k\ge 1}\left|\frac{X^N_k(0)}{N}-c_k(0)\right|\right){=}0
  \]
  hold. 
\item[(d)]
The initial state of the centered process, the random variables $W^N(0)$, $N{\ge}1$,  defined by relation~\eqref{Wn} satisfies the relation
\[
C_0{:=}\sup_N\E\left(\|W^N(0)\|_{L_2(w)}^2\right)<{+}\infty,
\]
and  the sequence $(W^N(0))$ converges in probability to an ${L_2(w)}$-valued random variable $W(0)$.
\end{itemize}
\end{assumption}
For a detailed  discussion of Assumption~(a), see Section~\ref{seccoeff} below. Assumption~(b) gives conditions on the moments of the initial state of the process. For example, by taking  $r_k{=}k^\beta$ with $\beta{>}1$, one can study the central limit problem in the Hilbert space $L_2(w)$ with weights $w_k{=}k^\alpha$ for some $\alpha{\in} (1,\beta)$. Note that this condition is more, but not much more,  demanding than the conservation of mass relation
\[
\sum_{k\ge 1}kX_k^N(t)/N\equiv 1,
\]
which is always satisfied. This assumption gives therefore a quite large class of initial distributions for which a central limit theorem holds.

\medskip
\paragraph{{\bf Some Notations.}}
For any $T{>}0$, let $\cal{D}_T{:=}D([0,T],{L_2(w)})$ be the space of c\`{a}dl\`{a}g functions on $[0,T]$ taking values in ${L_2(w)}$.
Since ${L_2(w)}$ is a separable and complete space, there exists a metric on $\cal{D}_T$, such that $\cal{D}_T$ is a separable and complete space.
See the chapter 3 in Billingsley~\cite{billing} for details.

Let $\cal{L}({L_2(w)})$ be the set of linear operators on ${L_2(w)}$ and the $\|\cdot\|_{\cal{L}(L_2(w))}$ be the associated norm for linear operators, \emph{i.e.} for any $f{\in}\cal{L}({L_2(w)})$,
\[\|f\|_{\cal{L}({L_2(w)})}=\sup_{\substack{z\in L_2(w)\\\|z\|_{L_2(w)}\leq 1}}\|f(z)\|_{L_2(w)}.\]

 For any $f,g\in \cal{L}({L_2(w)})$, $f{\circ} g$ denotes the composition. We call $A{\in} \cal{L}({L_2(w)})$ to be in a trace class if there exist a constant $C$, such that for all $z{\in} {L_2(w)}$, $\|A z\|_{L_2(w)}{\le} C\|z\|_{L_2(w)}$ and it has a finite trace, i.e.,
\[\textrm{Tr}A{:=}\sum_{n\ge 1}\left\langle A h_n,h_n\right\rangle_{L_2(w)}{<}\infty.\]

\section{SDE in Hilbert space}\label{limit}
In this section, we are going study the existence and uniqueness of  process
\begin{multline}\label{limitsde}
  W(t){=}W(0){+}\int_0^t\tau\bigg(\nabla s(c(u))\cdot W(u)\bigg)\diff u\\{+}\int_0^t\tau\bigg(\diag\left(\sqrt{s(c(u))}\right)\cdot\diff \beta (u)\bigg),
\end{multline}
when Assumptions~(a) and~(c) hold. 

Clearly, the process $(W(t))$ has the form of a stochastic differential equation in the infinite dimensional space $\R^{\N^+}$.
It has a stochastic integral
  \[
\cal{M}(t){:=}\int_0^t\tau\left(\diag\left(\sqrt{s(c(u))}\right)\cdot\diff \beta (u)\right)
\]
with respect to $\beta(t){=}\sum_{n\ge 1}\beta_n(t)\vec{e}_n$, which is a cylindrical Wiener process (c.f.  Yor~\cite{yor1974}) in the Hilbert space $L^2(\R^{\N^+})$.
Here $(c(t))$ is the unique solution of Becker-D\"oring equation~\eqref{bdeq} with initial state $c(0)$.
However, the mapping $\tau$ is unbounded in $L^2(\R^{\N^+})$.
Therefore, we consider the SDE~\eqref{limitsde} in the Hilbert space ${L_2(w)}$ (Definition~\ref{defH}).
In the following, we are going to show that 
\begin{enumerate}
\item
  for any $c{\in} \cal{X}_+$, the linear mapping $\tau{\circ} \nabla s(c): {L_2(w)}{\to} {L_2(w)}$  is bounded (and therefore Lipschitz);
  \item the stochastic process $(\cal{M}(t))$  is  well-defined   in ${L_2(w)}$ and is a martingale.
\end{enumerate}

Finally, by using the results in Section~7.1 of Da Prato and Zabczyk~\cite{Prato1992}, we give the proof of the existence and uniqueness of the solutions of the equation~\eqref{limitsde}.

 \begin{lemma}\label{tau}
    $\tau$ is a continuous linear mapping from ${L_2(w)}$ to ${L_2(w)}$: there exists a finite constant $\gamma_{\tau}(w)$, such that for any $z\in {L_2(w)}$,
    \[
\|\tau(z)\|_{L_2(w)}{\le} \gamma_{\tau}(w)\|z\|_{L_2(w)}.
    \]
\end{lemma}
\begin{proof}
  One only needs to verify that $\tau$ is bounded.
  For any $z\in {L_2(w)}$,
  \begin{equation*}
  \|\tau(z)\|_{L_2(w)}^2{=}R_1{+}R_2,
  \end{equation*}
  where
  \begin{align*}
    &R_1{=}w_1\left({-}\sum_{k{\ge} 1}(1{+}\ind{k{=}1})z_{2k{-}1}{+}\sum_{k{\ge} 2}(1{+}\ind{k{=}2})z_{2k{-}2}\right)^2,\\
    &R_2{=}\sum_{k{\ge} 2}{w_{k}}\left(z_{2k{-}3}{-}z_{2k{-}2}{-}z_{2k{-}1}{+}z_{2k}\right)^2.
  \end{align*}
  By using Cauchy-Schwarz inequality, we have
  \begin{equation*}
    R_1\le 4w_1\sum_{k{\ge} 1}\frac{1}{w_{k}}\|z\|_{L_2(w)}^2{=}4w_1\left\|\frac{1}{w}\right\|_{l_1}\|z\|_{L_2(w)}^2, 
      \end{equation*}
  and by using the increasing property of weights $w$
  \begin{multline*}
    R_2{\le} 4 \sum_{k{\ge} 2}w_{k}\left(z_{2k{-}3}^2{+}z_{2k{-}2}^2{+}z_{2k{-}1}^2{+}z_{2k}^2\right)\\
    {\le} 4\sum_{k{\ge} 2}\left(w_{2k{-}3}z_{2k{-}3}^2{+}w_{2k{-}2}z_{2k{-}2}^2
        {+}w_{2k{-}1}z_{2k{-}1}^2{+}w_{2k}z_{2k}^2\right){+}w_2z_1^2\\
        {\le }(8{+}w_2/w_1)\|z\|_{L_2(w)}^2.
    \end{multline*}

In conclusion, one has
$ \|\tau(z)\|_{L_2(w)}\le \gamma_{\tau}(w)\|z\|_{L_2(w)},$ for
\[
\gamma_{\tau}(w){=}\left(4w_1\left\|\frac{1}{w}\right\|_{l_1}{+}8{+}\frac{w_2}{w_1}\right)^{1/2}.
\]
\end{proof}

\begin{lemma}\label{drift}
  Under Assumption~(a), for any $c\in \cal{X}_+$ that satisfies
$\sum_{k\ge 1}{w_{k}}c_k{<}\infty,$ there exists a finite constant $\gamma(c,w)$, such that
  for any $x{\in} {L_2(w)}$,
\[
\|\nabla s(c)\cdot x\|_{L_2(w)}\le \gamma(c,w)\|x\|_{L_2(w)}.
\]
Moreover, under Assumption (a) and (c), if $(c(t))$ is the solution of the Becker-D\"oring equation~\eqref{bdeq} with initial state $c(0)$, then for any finite time $T$, one has
\[
\sup_{u\le T}\gamma(c(u),w){<}{+}\infty.
\]
\end{lemma}

\begin{proof}
  \begin{multline*}
    \|\nabla s(c)\cdot x\|_{L_2(w)}^2{=}\sum_{k{\ge} 1}w_{2k{-}1}a_k^2(x_1c_k{+}c_1x_k)^2{+}\sum_{k{\ge} 2}w_{2k{-}2} b_k^2x_k^2\\
    \le 2\gamma_{0}\Lambda^2\left(\sum_{k\ge 1}w_kc_k^2\right)x_1^2{+}2\gamma_0\Lambda^2\sum_{k{\ge} 1}{w_{k}}x_k^2{+}\gamma_0\Lambda^2\sum_{k{\ge} 2}{w_{k}}x_k^2\le \gamma(c,w)\|x\|_{L_2(w)}^2,
  \end{multline*}
  where for any $c\in\cal{X}_+$
  \[
\gamma(c,w){=}\gamma_0\Lambda^2\left(3{+}2\sum_{k\ge 1}w_kc_k\right).
\]
The Theorem~2.2 in Ball et al~\cite{BDeq1935} gives that $\sup_{u\le T}w_kc_k(u){<}\infty$ under Assumption~(a) and~(c).
\end{proof}

\subsection{Remark on the assumptions on the coefficients $(a_i)$ and $(b_i)$}\label{seccoeff}
We recall the main results for the existence and uniqueness of a solution of the deterministic Becker-D\"oring ODEs. 
\begin{enumerate}
\item From Theorem~2.2 in Ball et al.~\cite{Ball1986}, existence and uniqueness hold under the condition
  \[
a_i(g_{i+1}-g_i)=O(g_i)\qquad\textrm{  and }\qquad\sum_{i=1}^{+\infty} g_i c_i(0){<}{+}\infty,
\]
for a positive increasing sequence $g$ satisfying $\min_{i}(g_{i+1}-g_i){\ge }\delta{>}0$. For example, when $g_i{=}i^2$, \emph{i.e.}, $a_i{=}O(i)$ and the second moment of the initial state is bounded, the Becker-D\"oring equation is well-posed.
\item Conditions from Theorem~2.1 in Lauren\c{c}ot and Mischler~\cite{Laurencot2002},\\
  There exists some constant $K{>}0$ such that
  \[
  a_i{-}a_{i-1}{<}K \text{ and } b_i{-}b_{i-1}{<}K,
  \]
  for any $i{\in}\N$.
\end{enumerate}
It should be noted that  Lauren\c{c}ot and Mischler~\cite{Laurencot2002} have a stronger conditions on the coefficients but, contrary to Ball et al.~\cite{Ball1986}, not any on the initial state.

The known conditions necessary to get the convergence of the first order of the stochastic Becker-D\"oring model to the solution of the deterministic Becker-D\"oring ODEs are more demanding, they are given by Theorem~2 in Jeon~\cite{Jeon1998}, they are
\[
\lim_{i\to\infty}\frac{a_i}{i}{=}0\qquad\textrm{and}\qquad\lim_{i\to\infty}b_i{=}0.
\]
Note that a law of large numbers in this paper requires growth rates $(a_i)$ to be sublinear and  break rates $(b_i)$ to be vanishing.

For second order convergence, additional conditions seem however to be  necessary to establish the well-posedness of the limiting fluctuation process defined by SDE~\eqref{limitsde1}. 
First of all, for any $c{\in}\cal{X}_+$, the linear operator $\tau{\circ}\nabla s(c)$ does not seem to have monotonicity or symmetry properties that could give an alternative construction of the solution of SDE~\eqref{limitsde1}. Hence, boundedness properties of the linear operator $\tau{\circ}\nabla s(c)$ have to be used to get existence and uniqueness results of the solution of SDE~\eqref{limitsde1}.

For any $k{\ge}2$, it is easy to check that the $k$th coordinate of the vector  $\tau{\circ}\nabla s(c)\cdot h_{k}$ is ${-}(b_k{+}a_kc_1)/\sqrt{w_k}$, in particular the relation 
\[
\left\|\tau{\circ}\nabla s(c)\right\|_{\cal{L}(L_2(w))}\ge
\|\tau{\circ}\nabla s(c)\cdot h_{k}\|_{L_2(w)}\ge (b_k{+}a_kc_1)^2\qquad \textrm{for all }k,
\]
holds. Hence the operator $\tau{\circ}\nabla s(c)$ is unbounded in any weighted $L_2$ space if a boundedness property for the sequences $(a_i)$ and $(b_i)$ does not hold. By using similar arguments, we observe that the operator $\tau{\circ}\nabla s(c)$ is unbounded in the state spaces $l_1$ and $l_\infty$ as well if this property does not hold.

\begin{prop}\label{mg}
If Assumption (a)  and (c) hold, then the process $(\cal{M}(t))$ is a well-defined continuous, square-integrable, martingale.
\end{prop}
\begin{proof}
  By definition of the stochastic integral with respect to a cylindrical Wiener process (Chapter 4 in Da Prato and Zabczyk~\cite{Prato1992}), it is sufficient to verify that,
  \[
\int_0^T\sum_{n{\ge} 1}\|\tau{\circ} \diag\left(\sqrt{s(c(u))}\right)\cdot e_n\|_{L_2(w)}^2\diff u{<}\infty.
  \]
  By using Lemma \ref{tau},
  \begin{multline*}
  \|\tau{\circ} \diag\left(\sqrt{s(c(u))}\right)\cdot e_n\|_{L_2(w)}^2\\ \le \gamma_{\tau}(w)^2\|\diag\left(\sqrt{s(c(u))}\right)\cdot e_n\|_{L_2(w)}^2
              {=}{\gamma_{\tau}(w)^2}{w_{n}}s_n(c(u)).
\end{multline*}
  Therefore,
  \begin{align*}
    \int_0^T\sum_{n{\ge} 1}&\|\tau{\circ} \diag\left(\sqrt{s(c(u))}\right)\cdot e_n\|_{L_2(w)}^2\diff u\\
  &  \le {\gamma_{\tau}(w)^2}\int_0^T\sum_{k{\ge} 1}\left(w_{2k{-}1}a_kc_1(u)c_k(u){+}{w_{2k}}b_{k+1}c_{k+1}(u)\right)\diff u\\
 &\leq \gamma_0{\gamma_{\tau}(w)^2}\Lambda T\sup_{s\le T}\sum_{k{\ge} 1}{w_{k}}c_k(s){<}\infty.
  \end{align*}
  The last inequality is valid under Assumption (a) and (c).   See Ball et al.~\cite{BDeq1935} for details.
\end{proof}

\begin{theorem}
  For any measurable ${L_2(w)}$-valued random variable  $W(0)$, the  SDE~\eqref{limitsde} has a unique strong solution in ${L_2(w)}$.
  If, in addition, for $p{\ge} 1$, $\E\|W(0)\|_{L_2(w)}^{2p}{<}\infty$, then
  \[
\E\left(\sup_{0\le t\le T} \|W(s)\|_{L_2(w)}^{2p}\right){<}\infty.
  \]
\end{theorem}
\begin{proof}
  In Proposition \ref{mg}, we proved that the martingale part  $\cal{M}(t)$ of the SDE~\eqref{limitsde} is well-defined and continuous.
  By using Lemma \ref{tau} and Lemma \ref{drift}, the drift part is Lipschitz and linear on any $x\in C([0,T], {L_2(w)})$:
  \[
\|\tau\left(\nabla s(c(u))\cdot x(u)\right)\|_{L_2(w)}\le \gamma_\tau(w)\sup_{0\le u\le T}\gamma(c(u),w)\sup_{0\le u\le T}\|x(u)\|_{L_2(w)}.
  \]
 Therefore, by using the results in Section~7.1 of Da Prato and Zabczyk~\cite{Prato1992}, we can obtain the strongly existence, uniqueness, continuous and bounded results of SDE~\eqref{limitsde}.
\end{proof}
In Ball et al.~\cite{Ball1986} it is shown that, under the  assumptions
\begin{align}\label{as2}
  z_s^{{-}1}{:=}\limsup_{i\to\infty} R_i^{1/i}{<}\infty\qquad\textrm{and } \sup_{0\le z{<} z_s}\sum_{i{\ge} 1}iR_iz^i{>}1,
\end{align}
 where, for $k{\ge}1$,  $R_k{=}\prod_{i{=}2}^i(a_{i{-}1}/b_i)$,
then the equation $$\sum_{k{\ge} 1}kR_kz^k{=}1$$ has a unique solution $z{=}\widetilde{c}_1$. Moreover if $\widetilde{c}_k{=}R_k(\widetilde{c}_1)^i$, then $(\widetilde{c}_i)_{i{\ge} 1}$ 
is  the unique fixed point of Becker-D\"oring equations~(BD).

\begin{prop}[Fluctuations at equilibrium]\label{fe}
Under condition~\eqref{as2}, at equilibrium, the strong solution of SDE~\eqref{limitsde} can be represented as
  \begin{align}\label{sol}
\widetilde{W}(t){=}\cal{T}(t)\widetilde{W}(0){+}\int_0^t\cal{T}(t{-}s)\tau \bigg(\diff \cal{B}(s)\bigg),
  \end{align}
  where $(\cal{T}(t))$ is the  semi-group associated with linear operator $\tau{\circ} \nabla s(\widetilde{c})$ and $(\cal{B}(t))$ is a $\widetilde{Q}$-Wiener process in ${L_2(w)}$ where
  \[
\widetilde{Q}{=} \diag\left({w_{n}}s_n(\widetilde{c}),n{\ge} 1\right).
  \]
Moreover, the stochastic convolution part
\[
\widetilde{W}_{sc}(t){:=}\int_0^t\cal{T}(t-s)\tau\diff\cal{B}(s)
\]
is Gaussian, continuous in mean square. Its  previsible increasing processe is given by
\[
\left(\left\langle\widetilde{W}_{sc}(t)\right\rangle\right){=}\left(\int_0^t\cal{T}(r)\tau \widetilde{Q}\left(\cal{T}(r)\tau\right)^* \diff r\right).
\]
\end{prop}
\begin{proof}
At equilibrium, the SDE~\eqref{limitsde} becomes
\begin{align}\label{equsde}
   \widetilde{W}(t){=}\widetilde{W}(0){+}\int_0^t\tau\left(\nabla s(\widetilde{c})\cdot \widetilde{W}(u)\right)\diff u{+}\int_0^t\tau\left(\diag\left(\sqrt{s(\widetilde{c})}\right)\cdot\diff \beta (u)\right).
\end{align}
It is a linear equation with additive noise. The noise part can be expressed by
$\tau(\cal{B}(t))$ where
\[
\cal{B}(t){:=}\int_0^t \diag\left(\sqrt{s(\widetilde{c})}\right)\cdot\diff \beta (u){=}\sum_{n{\ge} 1}\sqrt{w_{n}s_n(\widetilde{c})}\beta_n(t)h_n.
\]
Recall that $(h_n)$ defined by~\eqref{defhn} is an orthonormal basis in ${L_2(w)}$.
By definition of $z_s$ and the assumption $\lim_{n\to\infty}w_n^{1/n}{=}1$, if $0\le z{<}z_s$, then 
\[
\sum_{n{\ge} 1}{w_{n}}R_nz^n{<}\infty
\]
holds. Let
  \[
\widetilde{Q}{=} \diag\left({w_{n}}s_n(\widetilde{c}),n{\ge} 1\right),
  \]
then
it is trace class, i.e.,
\[
\textrm{Tr}\widetilde{Q}{=}\sum_{n{\ge} 1}{w_{n}}s_n(\widetilde{c})\le 2\gamma_0\Lambda \sum_{n{\ge} 1}{w_{n}}\widetilde{c}_n{<}\infty.
\]
Therefore, $\cal{B}(t)$ is a $\widetilde{Q}$-Wiener process in ${L_2(w)}$ (c.f.
Section 4.1 in~\cite{Prato1992}).

Again, by using Lemma~\ref{tau} and Lemma~\ref{drift},
we can see that the linear operator $\tau{\circ} \nabla s(\widetilde{c})$ is bounded,
and therefore, the associated semi-group
\[
\cal{T}(t){:=}e^{\tau{\circ} \nabla s(\widetilde{c})t}
\]
is  uniformly continuous and satisfies that for any $z\in {L_2(w)}$,
\[
\|\cal{T}(t)z\|_{L_2(w)}\le e^{\gamma_\tau(w)\gamma(\widetilde{c},w) t}\|z\|_{L_2(w)}.
\]
Therefore, on any finite time interval $[0,T]$,
\begin{multline*}
  \int_0^t\textrm{Tr}\left(\cal{T}(r)\tau \widetilde{Q} \tau^*\cal{T}^*(r)\right)\diff r{=}\int_0^t\|\cal{T}(r)\tau\widetilde{Q}^{1/2}h_n\|_{L_2(w)}^2\diff r\\
  {\le} {\gamma_{\tau}(w)^2}Te^{2\gamma_{\tau}(w)\gamma(\widetilde{c},w) T}\textrm{Tr}{\widetilde{Q}}{<}\infty.
\end{multline*}

By applying the Theorem~5.4 in Da Prato and Zabczyk~\cite{Prato1992},
the process~\eqref{sol} is the unique weak solution of
SDE~\eqref{equsde}.
By Theorem~5.2 in~\cite{Prato1992}, the stochastic convolution part
$(\widetilde{W}_{sc}(t))$
is Gaussian, continuous in mean square, has a predictable version and for all $t\in [0,T]$,
\[
\left\langle\widetilde{W}_{sc}(t)\right\rangle{=}\int_0^t\cal{T}(r)\tau \widetilde{Q} \tau^*\cal{T}^*(r)\diff r.
\]
To show that the process defined by relation~\eqref{sol} is also a strong solution, we use Theorem~5.29 in~\cite{Prato1992}.
It is sufficient to check 
\[
\sum_{n{\ge} 1}\|\tau\nabla s(\widetilde{c})\tau \widetilde{Q}^{1/2}h_n\|_{L_2(w)}^2{<}\infty.
\]
Then,
stochastic integral $\tau{\circ}\nabla s(\widetilde{c})\cdot\widetilde{W}_{sc}(t)$ is a well defined continuous square-integrable martingale on $[0,T]$.
The process $\tau{\circ}\nabla s(\widetilde{c})\cdot\widetilde{W}(t)$ has square integrable trajectories. Therefore, $(\widetilde{W}(t))$ is a strong solution.
The proposition is proved. 
\end{proof}

\section{Convergence of the Fluctuation Processes}\label{convergence}
Recall that the fluctuation processes $(W^N(t))$ satisfies relation~\eqref{fluct},
\begin{multline*}
  W^N(t){=}W^N(0){+}\frac{1}{\sqrt{N}}\tau\left(D^N(t)\right){+}\frac{2a_1}{\sqrt{N}}\int_0^t\left(\frac{X_1^N(u)}{N}\right)e_1\diff u\\
  {+}\frac{1}{2}\int_0^t\tau\left(\nabla s\left(\frac{X^N(u)}{N}\right)\cdot W^N(u){+}\nabla s\left(c(u)\right)\cdot W^N(u)\right)\diff u.
\end{multline*}
The goal of this section is to prove that,  when $N$ is going to infinity,  the process $(W^N(t))$ is converging in distribution to  $(W(t))$, the solution of the SDE~\eqref{limitsde}.
We will   first prove some technical lemmas and the convergence of scaled process $(X^N(t)/N)$.
Then, we will prove the tightness and convergence of the local martingales $(D^N(t))$ and, with the help of these results, we will get the tightness of $(W^N(t))$ and identify the limit.

\begin{lemma}
Under Assumptions~(a) and~(b) then,   for any $T{>}0$, one has
  \begin{align}
    \zeta_T&{:=}\sup_N\sup_{t{\le} T}\E\left(\sum_{k{\ge} 1}\frac{r_kX_k^N(t)}{N}\right){<}{+}\infty,\label{zetaT}\\
    \kappa_T&{:=}\sup_N\E\left(\sup_{t{\le} T}\sum_{k{\ge} 1}\frac{\left|\cal{M}_k^N(t)\right|}{\sqrt{N}}\right){<}{+}\infty,\label{kappaT}
  \end{align}
    where $\cal{M}^N(t){=}\tau(D^N(t))$.
\end{lemma}
\begin{proof}
  By using the SDE~\eqref{bdsde},
  we get that, for any $N$,
  \begin{multline*}
    \E\left(\sum_{k{\ge} 1}\frac{r_kX_k^N(t)}{N}\right){=}
    \E\left(\sum_{k{\ge} 1}\frac{r_k X_k^N(0)}{N}\right)\\
           {+}\E\int_0^tr_{2}\left(a_1\frac{X_1^N(u)(X_1^N(u){-}1)}{N^2}{-}b_2\frac{X_2^N(u)}{N}\right)\diff u\\
    {+}r_{1}
    \E\int_0^t\bigg({-}\sum_{i{\ge} 1}(1{+}\ind{i{=}1})a_{i}\frac{X_1^N(u)(X_i^N(u){-}\ind{i{=}1})}{N^2}{+}\sum_{i{\ge} 2}(1{+}\ind{i{=}2})b_i\frac{X^N_i(u)}{N}\bigg)\diff u\\
    {+}\E\int_0^t\bigg(\frac{X_1^N(u)}{N}\sum_{i{\ge} 2}(r_{i+1}{{-}}r_{i})a_i\frac{X_i^N(u)}{N}{+}\sum_{i{\ge} 3}(r_{i{-}1}{-}r_{i})b_i\frac{X^N_i(u)}{N}\bigg)\diff u\\
    {\le}  \E\left(\sum_{k{\ge} 1}\frac{r_kX_k^N(0)}{N}\right)  {+} \left(r_{2}{+}r_{1}\right)\Lambda t{+}\gamma_r\Lambda \E\int_0^t\sum_{i{\ge} 1}r_{i}\frac{X_i^N(u)}{N}\diff u,
  \end{multline*}
  where $\Lambda$ is defined in Assumption~(a). 
  We apply  Gronwall's inequality,
  then, for any $t{\le} T$,
  \begin{align*}
    \E\left(\sum_{k{\ge} 1}\frac{r_kX_k^N(t)}{N}\right){\le} e^{\gamma_r\Lambda T}\left(\E\left(\sum_{k{\ge} 1}\frac{r_k X_k^N(0)}{N}\right){+}(r_{2}+r_{1})\Lambda T\right).
  \end{align*}
  Therefore, there exists a constant $\zeta_T$, such that   the relation~\eqref{zetaT} holds.
  
  For the other inequality, we first see
  \begin{align*}
 &\E\left(\sup_{t{\le} T}\sum_{k{\ge} 1}\frac{\left|\cal{M}_k^N(t)\right|}{\sqrt{N}}\right)\\ &{\le}
    \sum_{k{\ge} 1}\frac{1}{{w_k}}{+}\E\left(\sup_{t{\le} T}\sum_{k{\ge} 1}\frac{\left|\cal{M}_k^N(t)\right|}{\sqrt{N}}\ind{\frac{\left|\cal{M}_k^N(t)\right|}{\sqrt{N}}{\ge} \frac{1}{{w_k}}}\right)\\
    &{\le}
   2{+}\E\left(\sup_{t{\le} T}\sum_{k{\ge} 1}\frac{{w_k}\cal{M}_k^N(t)^2}{N}\right)
    {\le}
   2{+}\E\left(\sum_{k{\ge} 1}\frac{{w_k}\langle\cal{M}_k^N(T)\rangle}{N}\right).
  \end{align*}
 By using the expression of the increasing process~\eqref{inc} and the definition of $\tau$, we get
  \begin{multline*}
\sum_{k{\ge} 2}\frac{{w_k}}{N}\left\langle\cal{M}_k^N(T)\right\rangle{\le}
\Lambda\int_0^T\sum_{k{\ge} 1}(w_{k{+}1}\ind{k{\ge} 2}{+}{w_k})\frac{X_{k}^N(u)}{N}\diff u\\
{+}\Lambda\int_0^T\sum_{k{\ge} 2}({w_k}{+}w_{k{-}1}\ind{k{\ge} 3})\frac{X_{k}^N(u)}{N}\diff u.
  \end{multline*}
  It implies
  \begin{align*}
\sup_N\E\left(\sum_{k{\ge} 2}\frac{{w_k}}{N}\left\langle\cal{M}_k^N(T)\right\rangle\right){\le}
w_1\Lambda T{+}(3{+}\gamma_0)\Lambda T\sup_N\sup_{t{\le} T}\E\sum_{k{\ge} 2}\frac{{w_k}}{N}X_{k}^N(t)\\
{\le}
w_1\Lambda T{+}(3{+}\gamma_0)\Lambda T\zeta_T.
  \end{align*}
  For $k{=}1$,
  \[
\left\langle\frac{\cal{M}_1^N(T)}{N}\right\rangle{\le} 2\int_0^T\sum_{i{\ge} 1}\left(a_i\frac{X_1^N(u)X_i^N(u)}{N^2}{+}b_{i{+}1}\frac{X_{i{+}1}^N(u)}{N}\right)\diff u{\le} 4\Lambda T.
  \]
  Therefore,  there exists a positive constant $\kappa_T$, such that
\[
\sup_N\E\left(\sup_{t{\le} T}\sum_{k{\ge} 1}\frac{\left|\cal{M}_k^N(t)\right|}{\sqrt{N}}\right){\le} \kappa_T
\]
holds. 

\end{proof}


Now we prove the law of large numbers under Assumptions~(a),~(b) and~(c)
in the $L_1$-norm. This result will be used in the proof of Theorem~\ref{fclt}.
We should note that we could not directly apply the Theorem 2 of Jeon~\cite{Jeon1998} since it requires the conditions $\lim_{i\to\infty} a_i/i{=}0$ and $\lim_{i{\to} \infty}b_i{=}0$.

\begin{theorem}\label{lln}
  {\normalfont \bf (Law of large numbers)}
  Under Assumptions~{\normalfont(a)}, {\normalfont(b)} and {\normalfont(c)},
for any $T{>}0$, one has
  \[
\lim_{N\to\infty}\E\left(\sup_{t{\le} T}\sum_{k{\ge} 1}\left|\frac{X^N_k(t)}{N}{-}c_k(t)\right|\right){=}0,
\]
where $(c(t))$ is the unique solution of Becker-D\"oring equation~{\normalfont(BD)} with initial state $c(0){\in}\cal{X}_{+}$.
\end{theorem}

\begin{proof}
From SDE~\eqref{bdsde}, we have that
  \begin{multline*}
    \frac{X^N(t)}{N}{-}c(t){=}
    \frac{X^N(0)}{N}{-}c(0)
    {+}\int_0^t\tau\left(s\left(   \frac{X^N(u)}{N}\right){-}s\left(   c(u)\right)\right)\diff u\\
    {+}\frac{2a_1}{N}\int_0^t\left(\frac{X_1^N(u)}{N}\right)e_1\diff u{+}\frac{\cal{M}^N(t)}{N}.
  \end{multline*}
  With a simple calculation, we get the relation
  \begin{multline*}
\sum_{k{\ge} 1}\left|\frac{X^N_k(t)}{N}{-}c_k(t)\right|{\le} 
\sum_{k{\ge} 1}\left|\frac{X^N_k(0)}{N}{-}c_k(0)\right|{+}\frac{2\Lambda t}{N}
\\{+}\int_0^t8\Lambda \sum_{k{\ge} 1}\left|\frac{X^N_k(u)}{N}{-}c_k(u)\right|\diff u
{+}\sum_{k{\ge} 1}\left|\frac{\cal{M}_k^N(t)}{N}\right|,
  \end{multline*}
and,   by Gronwall's inequality and relation~\eqref{kappaT}, we have
  \[
\lim_{N\to\infty}\E\sup_{t\le T}\sum_{k{\ge} 1}\left|\frac{X^N_k(t)}{N}{-}c_k(t)\right|{=}0.
\]
\end{proof}

\begin{prop}\label{ppd}
Under Assumptions~(a) and (b),   the sequence of local martingales 
  \[
\left(\frac{1}{\sqrt{N}}D^N(t)\right){=}\left(\frac{1}{\sqrt{N}}D_i^N(t)\right)_{i{\ge} 1}
\]
 is tight for the convergence in distribution in $\cal{D}_T$.
\end{prop}
\begin{proof}

  Recall that for all $i{>}N$, $(D_i^N(t))\equiv 0$ and $(X^N_i(t))\equiv 0$.
  By using Jensen's and Doob's inequalities, for any fixed $T$ and $N$, we get
  \begin{multline*}
   \E \left(\sup_{t{\le}T}\left\|\frac{1}{\sqrt{N}}D^N(t)\right\|_{L_2(w)}\right)\le
    \sqrt{\E\left(\sup_{t{\le}T}\left\|\frac{1}{\sqrt{N}}D^N(t)\right\|_{L_2(w)}^2\right)}\\
    {=}\sqrt{\E\left(\sup_{t{\le}T}\sum_{i{=} 1}^N {w_i} \frac{1}{N}D^N_i(t)^2\right)}
    \le
    \sqrt{\sum_{i{=} 1}^N {w_i}\frac{1}{N}\E\left(\langle D^N_i(T)\rangle\right)}.
  \end{multline*}
From the expression of the increasing processes~\eqref{inc}, we have
  \[
\frac{1}{N}\sum_{k{\ge} 1}w_{2k{-}1}\E\left(\langle D^N_{2k{-}1}(T)\rangle\right){\le}\gamma_0\Lambda\int_0^T\E\left(\sum_{k{\ge} 1}\frac{{w_k}X_k^N(u)}{N}\right)\diff u,
\]
and
  \[
\frac{1}{N}\sum_{k{\ge} 2}w_{2k{-}2}\E\left(\langle D^N_{2k{-}2}(T)\rangle\right){\le}\gamma_0\Lambda\int_0^T\E\left(\sum_{k{\ge} 2}\frac{{w_k}X_k^N(u)}{N}\right)\diff u.
\]
Therefore, by using inequality~\eqref{zetaT}, one gets
\[
 \sup_N\E\sup_{t{\le}T}\left\|\frac{1}{\sqrt{N}}D^N(t)\right\|_{L_2(w)}{\le}\sqrt{2\gamma_0\Lambda T\sup_{N}\sup_{t{\le}T}\E\left(\sum_{k{\ge} 1}\frac{{w_k}X_k^N(t)}{N}\right)}{<}\infty,
 \]
 it gives that 
\begin{multline*}
 \lim_{a\to\infty}\limsup_{N\to\infty}\P\left(\sup_{t{\le}T}\left\|\frac{1}{\sqrt{N}}D^N(t)\right\|_{L_2(w)}{\ge} a\right)\\{\le}\lim_{a\to\infty}\frac{1}{a}\sup_N\E\left(\sup_{t{\le}T}\left\|\frac{1}{\sqrt{N}}D^N(t)\right\|_{L_2(w)}\right){=}0.
\end{multline*}
 Since the Skorohod distance is weaker than the uniform distance in $\cal{D}_T$, we estimate the modulus of continuity with uniform distance to prove the tightness in the Skorohod space. For any $\varepsilon{>}0$, $\delta{>}0$ and any $L{\in}\N$, we have 
\begin{align*}
  \P\left(\sup_{t,s{\le}T,|t{-}s|{<}\delta}\left\|\frac{1}{\sqrt{N}}D^N(t){-}\frac{1}{\sqrt{N}}D^N(s)\right\|_{L_2(w)}{\ge} \varepsilon\right)
  \le\P\left(
  \sup_{t{\le}T}\sum_{i{>}L}r_{i}\frac{1}{N}D_i^N(t)^2{\ge} \frac{\varepsilon^2 r_{L}}{4w_L}
  \right)\\
  {+}\P\left(\sup_{t,s{\le}T,|t{-}s|{<}\delta}\sum_{i{=}1}^L{w_{i}}\left(\frac{1}{\sqrt{N}}D^N_i(t){-}\frac{1}{\sqrt{N}}D^N_i(s)\right)^2{\ge} \frac{\varepsilon^2}{2}\right).
\end{align*}
By using the inequality~\eqref{zetaT},
\begin{align*}
\P\left(
  \sup_{t{\le}T}\sum_{i{>} L}r_{i}\frac{1}{N}D_i^N(t)^2{\ge} \frac{\varepsilon^2 r_{L}}{4w_L}
  \right)\le
  \frac{4w_L}{\varepsilon^2 r_{L} }\E\left(
  \sum_{i{>} L}r_{i}\frac{1}{N}\langle D_i^N(T)\rangle
  \right)\\
  {\le}
  \frac{8\gamma_r w_L\Lambda }{\varepsilon^2 r_{L}}T\sup_{N}\sup_{t{\le}T}\E\left(
  \sum_{i{\ge} 1}\frac{r_{i}X_i^N(t)}{N}\right)\le\frac{8\gamma_r w_L\Lambda }{\varepsilon^2 r_{L}}T\zeta_T.
\end{align*}
By using the Assumption (b), $\lim_{k\to\infty}w_k/r_k=0$, for any constant $\eta{>}0$, there exist a constant $L$, such that
\[
\P\left(
  \sup_{t{\le}T}\sum_{i{\ge} L}r_{i}\frac{1}{N}D_i^N(t)^2{\ge} \frac{\varepsilon^2 r_L}{4w_L}
  \right){\le}\frac{\eta}{2}.
  \]
  The processes $(D_i^N(t)/\sqrt{N})$ $i{=}1$,\dots, $L$, live in finite dimensional space and each of them has an increasing process that is uniformly continuous almost surely.
Therefore, for $N$ large enough, we have
  \[
\P\left(\sup_{t,s{\le}T,|t{-}s|{<}\delta}\sum_{i{=}1}^L{w_i}\left(\frac{1}{\sqrt{N}}D^N_i(t){-}\frac{1}{\sqrt{N}}D^N_i(s)\right)^2{\ge} \frac{\varepsilon^2}{2}\right){\le}\frac{\eta}{2},
\]
consequently, by using Theorem 13.2 in Billingsley~\cite{billing}, the sequence of processes $(D^N(t)/\sqrt{N})$ is tight in $\cal{D}_T$.

\end{proof}

\begin{prop}\label{limitd}
  For the convergence in distribution of random process, uniformly on compact sets, 
  \[
\lim_{N\to{+}\infty} \left(\frac{1}{\sqrt{N}}D^N(t)\right){=}\left(D(t)\right){:=}\left(\int_0^t  \diag\left(\sqrt{s(c(u))}\right)\cdot\diff \beta(u)\right).
\]
\end{prop}
\begin{proof}

From the previous proposition, we have that, for $T{>}0$, the sequence $(D^N(t)/\sqrt{N})$ is tight in $D([0,T],{L_2(w)})$. Let $D'(t)$ be a possible limit.
For any $d{\in}\N^*$, let $\cal{P}_d$ be the projection from $\R^{\N^{+}}$ to $\R^d$, i.e., for any $z{\in}\R^{\N^{+}}$, $\cal{P}_d(z){=}(z_1,\dots,z_d).$ It is easy to check that $(\cal{P}_d(D'(t)))$ is a limit of a subsequence of
$(\cal{P}_d(D^N(t)/\sqrt{N}))$ for the weak convergence in probability in the $L_2$-norm.
By using Theorem~1.4 page 339 of Ethier and Kurtz~\cite{EAK}, we know  that for any $d{\in}\N^*$, the equality
$\left(\cal{P}_d(D'(t))\right){=}\left(\cal{P}_d(D(t))\right)$ holds in distribution. Hence, by Kolmogorov's theorem,  we have the equality in distribution $\left(D'(t)\right){=}\left(D(t)\right)$.
\end{proof}

 We can now state our main result. 
\begin{theorem}[Functional Central Limit Theorem]\label{fclt}
  Under Assumptions {\normalfont(a)--(d)}, then  the fluctuation process $(W^N(t))$ defined by relation~\eqref{Wn} converges  in distribution to the ${L_2(w)}$-valued process $(W(t))$, the unique strong solution of the SDE~\eqref{limitsde}.
\end{theorem}
\begin{proof}
 From relation~\eqref{fluct}, process $(W^N(t))$ satisfies
\begin{multline*}
    W^N(t){=}W^N(0){+}\frac{1}{\sqrt{N}}\tau\left(D^N(t)\right)
    {+}\int_0^t\tau\circ\nabla s\left(c(u)\right)\cdot W^N(u)\diff u\\
    {+}\frac{2a_1}{\sqrt{N}}\int_0^t\left(\frac{X_1^N(u)}{N}\right)e_1\diff u{+}\int_0^t\tau\left(\Delta^N(u)\right)\diff u
\end{multline*}
  where
  \begin{align*}
    \Delta^N(u){=}\frac{1}{2}\left(\nabla s\left(\frac{X^N(u)}{N}\right){+}\nabla s\left(c(u)\right)\right)\cdot W^N(u){-}\nabla s\left(c(u)\right)\cdot W^N(u).
  \end{align*}
  By direct calculations, we have that for $k{\ge} 1$
  \begin{align*}
    \Delta_{2k{-}1}^N(u){=}&\left(1{+}\ind{k{=}1}\right)a_k\left(\frac{X^N_1(u)}{N}{-}c_1(u)\right)W^N_k(u),\\
    \Delta_{2k}^N(u){=}&0,
  \end{align*}
and then
  \begin{align*}
\|\Delta^N(u)\|_{L_2(w)}{\le}2\sqrt{\gamma_0}\Lambda\left|\frac{X^N_1(u)}{N}{-}c_1(u)\right|\|W^N(u)\|_{L_2(w)}{\le}2\sqrt{\gamma_0}\Lambda\|W^N(u)\|_{L_2(w)}.
  \end{align*}
Let $\Gamma=\gamma_\tau(w)(\sup_{u\le T}\gamma(c(u),w){+}2\sqrt{\gamma_0}\Lambda)$, by using Lemma~\ref{tau} and Lemma~\ref{drift}, for all $t{\le}T$, one has
  \begin{multline*}
    \|W^N(t)\|_{L_2(w)}{\le} \|W^N(0)\|_{L_2(w)}{+}\gamma_{\tau}(w)\frac{1}{\sqrt{N}}\|D^N(t)\|_{L_2(w)}{+}\frac{2\sqrt{w_1}\Lambda t}{\sqrt{N}}\\
    {+}\int_0^t\Gamma\|W^N(u)\|_{L_2(w)}\diff u.
  \end{multline*}
  Thanks to Gronwall's lemma, for any $t{\le}T$, we have
    \begin{align*}
\|W^N(t)\|_{L_2(w)}{\le}e^{\Gamma T}\left(\|W^N(0)\|_{L_2(w)}{+}\gamma_\tau(w)\frac{1}{\sqrt{N}}\|D^N(t)\|_{L_2(w)}{+}\frac{2\sqrt{w_1}\Lambda T}{\sqrt{N}}\right),
    \end{align*}
    which implies the relation
  \begin{multline*}
    \sup_N\E\sup_{t{\le}T}\|W^N(t)\|_{L_2(w)}^2\\
        {\le} 3e^{2\Gamma T}\left(\sup_N\E\|W^N(0)\|_{L_2(w)}^2{+}\gamma_\tau(w)^2\sup_N\E\sup_{t{\le}T}\frac{1}{{N}}\|D^N(t)\|_{L_2(w)}^2{+}\frac{4w_1\Lambda^2T^2}{N}\right).
  \end{multline*}
In the proof of Proposition~\ref{ppd}, we have shown that
\[\sup_N\E\sup_{t{\le}T}\frac{1}{{N}}\|D^N(t)\|_{L_2(w)}^2{<}\infty,\]
therefore, if
\[
\sup_N\E\|W^N(0)\|_{L_2(w)}^2{<}C_0,
\]
then there exists a finite constant $C_T$ such that
\[
\sup_N\E\sup_{t{\le}T}\|W^N(t)\|_{L_2(w)}^2{<}C_T.
\]
  For the tightness of $(W^N(t))$, for any $\varepsilon{>}0$, $\eta{>}0$, by using Lemma~\ref{tau} and Lemma~\ref{drift}, we get that
  \begin{multline*}
    \mathbb{P}\left(\sup_{t,s{\le}T,|t{-}s|{<}\delta}\|W^N(t){-}W^N(s)\|_{L_2(w)}{\ge} \varepsilon\right){\le}\P\left(\Gamma\delta\sup_{u{\le}T}\|W^N(u)\|_{L_2(w)}{\ge} \frac{\varepsilon}{2}\right)\\
 {+} \mathbb{P}\left(\sup_{t,s{\le}T,|t{-}s|{<}\delta}\left\|\frac{D^N(t)}{\sqrt{N}}{-}\frac{D^N(s)}{\sqrt{N}}\right\|_{L_2(w)}{\ge} \frac{\varepsilon}{4\gamma_\tau(w)}\right){+}\ind{\frac{2\sqrt{w_1}\Lambda\delta}{\sqrt{N}}{>}\frac{\varepsilon}{4}}
  \end{multline*}
  holds. Choose $\delta_1{>}0$ such that
  \[
\delta_1{<}\sqrt{\frac{\eta}{2C_T}}\frac{\varepsilon}{2\Gamma},
\]
then, for  $\delta{\in}(0,\delta_1)$ and  $N{{\ge}}1$,
\[
\P\left(\Gamma\delta\sup_{u{\le}T}\|W^N(u)\|_{L_2(w)}{\ge} \frac{\varepsilon}{2}\right)
       {\le}\left(\frac{2\Gamma\delta}{\varepsilon}\right)^2\E\sup_{u{\le}T}\|W^N(u)\|_{L_2(w)}^2{<}\frac{\eta}{2}.
\]
According to the proof of Proposition~\ref{ppd}, there exist  $\delta_2{>}0$ and $N_0$, such that, for  $\delta{\in}(0,\delta_2)$ and $N{>}N_0$, the relation
\[
\mathbb{P}\left(\sup_{t,s{\le}T,|t{-}s|{<}\delta}\left\|\frac{D^N(t)}{\sqrt{N}}{-}\frac{D^N(s)}{\sqrt{N}}\right\|_{L_2(w)}{\ge} \frac{\varepsilon}{4\gamma_\tau(w)}\right){\le}\frac{\eta}{2}
\]
holds. In conclusion, for any $\delta{<}\delta_1{\land}\delta_2{\land}(\varepsilon/(8\Lambda \sqrt{w_1}))$ and $N{>}N_0$,
\begin{align*}
    \mathbb{P}\left(\sup_{t,s{\le}T,|t{-}s|{<}\delta}\|W^N(t){-}W^N(s)\|_{L_2(w)}{\ge} \varepsilon\right){\le}\eta.
  \end{align*}
  It is then easy to check that
  \[
\lim_{a\to\infty}\limsup_{N\to\infty}\P\left(\sup_{t{\le}T}\|W^N(t)\|_{L_2(w)}{\ge} a\right){=}0.
\]
Therefore, the process $(W^N(t))$ is tight in $D([0,T],{L_2(w)})$.
To identify the limit, note that
 \begin{multline*}
    \E\sup_{t{\le}T}\left\|\int_0^t\tau\left(\Delta^N(u)\right)\diff u\right\|_{L_2(w)}\le\gamma_\tau(w)\E\int_0^T\|\Delta^N(u)\|_{L_2(w)}\diff u\\
    {\le}2\sqrt{\gamma_0}\gamma_\tau(w)T\Lambda \E\left(\sup_{u{\le}T}\left|\frac{X^N_1(u)}{N}{-}c_1(u)\right|\|W^N(u)\|_{L_2(w)}\right)\\
    \le2\sqrt{\gamma_0}\gamma_\tau(w)T\Lambda\left(2\E  \sup_{u{\le}T}\left|\frac{X^N_1(u)}{N}{-}c_1(u)\right|\right)^{1/2}
    \left(\E  \sup_{u{\le}T}\|W^N(u)\|_{L_2(w)}^2\right)^{1/2}.
  \end{multline*}
  By using Theorem~\ref{lln}, this term is vanishing as $N$ goes to infinity.
From Proposition~\ref{limitd}, one conclude that any limit of $(W^N(t))$ satisfies SDE~\eqref{limitsde}. The theorem is proved. 
\end{proof}

\subsection*{Acknowledgments}
The author thanks two anonymous referees for their comments and references that have been very helpful to clarify several aspects of this paper. 
The author also would like to express her gratitude to her advisor Professor Philippe Robert for valuable discussions.
The author's work is supported by the grant from Fondation Sciences Math\'ematiques de Paris (FSMP),
overseen by the French National Research Agency (ANR) as part of the ``Investissements d'Avenir'' program (reference: ANR-10-LABX-0098).
\bibliographystyle{amsplain}
\bibliography{ref}
\end{document}